\documentclass[12pt]{amsart}
\usepackage{times}
\usepackage{amsfonts}
\usepackage{amsthm}
\usepackage{amsmath}
\advance \topmargin by -\headheight
\advance \topmargin by -\headsep
\evensidemargin \oddsidemargin
\newtheorem{theorem}{Theorem}[section]
\newtheorem{lemma}[theorem]{Lemma}
\newtheorem{corollary}[theorem]{Corollary}

\begin{document}

\title{Some Free Entropy Dimension Inequalities for Subfactors}

\author{Kenley Jung}

\address{Department of Mathematics, University of California,
Los Angeles, CA 90095-1555, USA}

\email{kjung@math.ucla.edu}
\thanks{Research supported by an NSF Graduate Fellowship and an NSF Postdoctoral Fellowship}

\begin{abstract} Suppose $N \subset M$ is an inclusion of $II_1$-factors of finite index.  If $N$ can be generated by a finite set of elements, then there exist finite generating sets $X$ for $N$ and $Y$ for $M$ such that $\delta_0(X) \geq \delta_0(Y),$ where $\delta_0$ denotes Voiculescu's microstates (modified) free entropy dimension.  Moreover, given $\epsilon>0$ one has 

\[ \delta_0(F) \geq \delta_0(G) \geq ([M:N]^{-2} - \epsilon) \cdot (\delta_0(F) -1)+1 - \epsilon \] 

\vspace{.1in}

\noindent for certain generating sets $F$ for $N$ and $G$ for $M.$     
\end{abstract}

\maketitle

Connections between free probability and subfactor theory have existed for
at least a decade.  In [7] Radulescu studied finite index subfactors of
the free group factors.  Stefan applied microstate techniques to study
finite index subfactors of the interpolated free group factors.  He showed
in [9] that finite index subfactors of free group factors have certain
irreducibility properties like primality (see also [6] where Ozawa uses
$C^*$ techniques to generalize these primality results).  More
recently Shlyakhtenko and Ueda have used free probability techniques to
construct irreducible subfactors of $L(F_{\infty})$ ([8]).

We are interested in the problem of computing the microstates free entropy dimension of generating sets of a finite inclusion of subfactors.  Suppose $N \subset M$ is a finite inclusion of $II_1$-factors and $X$ and $Y$ are finite sets of generators for $N$ and $M,$ respectively.  One would expect the following:

\[ \delta_0(X) = [M:N] (\delta_0(Y) -1) +1.\]  

\noindent This formula is a straighforward analogue of Shreier's theorem on the free groups: if $H \subset G$ is an inclusion of groups of finite index and $G$ is free group on $n$ generators, then $H$ is a free group on $[G:H](n-1) +1$  generators.

In this remark we find generating sets $X$ for $N$ and $Y$ for $M$ such that $\delta_0(X) \geq \delta_0(Y)$. This inequality is clearly weaker than the equality previously stated since $[M:N] \geq 1.$   The argument makes use of the basic construction in [3] and the hyperfinite inequality in [5].  Using this monotonicity result, duality in the basic construction, and a computation of the free entropy dimension for amplifications, it follows that for $\epsilon >0$ there exist generating sets $X$ for $N$ and $Y$ for $M$ such that

\[ \delta_0(X) \geq \delta_0(Y) \geq ([M:N]^{-2} -\epsilon) \cdot (\delta_0(X) -1) +1 -\epsilon.\]  

\noindent The $\epsilon$'s can be dropped provided that $[M:N] \in \mathbb Q.$  As a corollary it follows that if $L(F_r)$ is the interpolated free group with index $r$ and $M$ is a $II_1$-factor such that for some $\epsilon >0$ and any finite set of generators $X$ for $M$ $\delta_0(X) <r- \epsilon,$ then $L(F_r)$ cannot contain $M$ as a subfactor of finite index.

Section 1 is a brief list of preliminaries and assumptions.  Section 2 will prove the monotonicity result.  Section 3 deals with the free entropy dimension of amplifications and Section 4 will present the final inequality.  Section 5 consists of applications of the preceding material.

\section{Preliminaries and Assumptions}

Throughout suppose $M$ is a $II_1$-factor with tracial state $\varphi$ and $N \subset M$ is a subfactor.  We will always assume that $M$ embeds into the ultraproduct of the hyperfinite $II_1$-factor (this guarantees the existence of microstates for any set of generators for $M$). For any $n \in \mathbb N, |\cdot|_2$ denotes
the norm on $(M_k(\mathbb C))^n$ given by $|(x_1,\ldots, x_n)|_2
= ( \sum_{j=1}^n tr_k(x_j^2))^{\frac{1}{2}}$ where $tr_k$ is the
tracial state on the $k\times k$ complex matrices.  We will maintain the notation
introduced in [5] and [10].  If $F = \{a_1, \ldots, a_n\}$ is a
finite set of elements in $M,$ we abreviate
$\Gamma_R(a_1, \ldots, a_n;m,k,\gamma)$ by $\Gamma_R(F;m,k,\gamma)$
and in a similar way we write the associated microstate sets and
quantities $\delta_0(F), \mathbb
P_{\epsilon}(F),$ and $\mathbb K_{\epsilon}(F).$  We point out that $F$ does not necessarily consist of selfadjoint elements and thus, the microstate spaces are tacitly assumed to consist of matricial approximants \emph{in $*$-moments}.  This will cause no problems, since all these concepts make sense with non-selfadjoint elements (as has been observed by many before), and moreoever, clearly coincide with the selfadjoint quantities.   Finally, if $G = \{b_1,
\ldots, b_p\}$ is another finite set of elements in $M,$
then we denote by $\Gamma_R(F \cup G;m,k, \gamma)$ the set
$\Gamma_R(a_1, \ldots, a_n, b_1, \ldots, b_p;m,k,\gamma)$ and write
all the associated microstate quantities $\delta_0(F \cup G), \mathbb
P_{\epsilon}(F \cup G),$ and $\mathbb K_{\epsilon}(F \cup G)$ with
respect to $\Gamma_R(F \cup G;m,k,\gamma).$ Finally,
$\Gamma(F;m,k,\gamma)$ will denote the set of all $k\times k$
$*$-microstates (no restrictions on the operator norms)  with degree of
approximation $(m,\gamma).$

\section{A Monotic Property of $\delta_0$}

Recall the basic construction and some of its properties introduced in [3].  On $L^2(M)$ consider the orthogonal projection $e_N$ onto $L^2(N) \subset L^2(M).$  Denote by $M_1$ the von Neumann subalgebra of $B(L^2(M))$ generated by $M$ and $e_N.$  $M_1$ is called the basic construction for $N \subset M.$  $M_1$ turns out to be a type $II$-factor and when $[M:N] < \infty,$ then $M_1$ is a type $II_1$-factor.  

By [3] any inclusion of $II_1$-factors of finite index arise from the basic construction.  By this I mean that there exists a $II_1$-subfactor $N_1$ of $N$ of finite index such that the basic construction $P$ for $N_1 \subset N$ is isomorphic to $M$ in such a way that the two copies of $N$ (one in $P$ and one in $M$) are canonically identified.  The crucial point is the fact that the projection $f$ onto $L^2(N_1) \subset N$ together with a set of generators for $N$ is a generating set for $M$ (identified with $P$) in such a way that this set of generators for $M$ has an element (namely, $f$) which commutes with any element of $N_1$.

With these preliminaries out of the way we come to the main point of this section which more or less says that if $N$ can be finitely generated, then there exist finite generating sets $Y$ for $N$ and $Z$ for $M$ such that $\delta_0(Y) \geq \delta_0(Z).$  In order to deal with some technical issues in Section 4, I will prove below a bloated version of this simple statement.

\begin{lemma} If $X$ is a finite set of generators for $N,$ then there exist finite sets of selfadjoints $F$ and $G$ such each generates a copy of the hyperfinite $II_1$-factor, $F^{\prime \prime} \subset G^{\prime \prime},$ $X \cup F$ generates $N,$ $X \cup G$ generates $M,$ $f_0 \in G,$ and 

\[\delta_0(X \cup F) \geq \delta_0(X \cup G).\]

\end{lemma}

\begin{proof}   Suppose $\mathcal R$ is a copy of the hyperfinite $II_1$-factor in $N_1$ and pick a finite set of selfadjoint generators $F$ for $\mathcal R.$  $f$ commutes with $\mathcal R.$  Hence by the computations in [4] and the hyperfinite inequality for $\delta_0$ in [5]

\[ \delta_0(X \cup F \cup f) \leq \delta_0(X \cup F) + \delta_0(F \cup \{f\}) - \delta_0(F) = \delta_0(X \cup F) + 1 -1 = \delta_0(X \cup F).\]   

\noindent There exists a von Neumann algebra $\mathcal R_1 \subset M$ isomorphic to the hyperfinite $II_1$-factor such that $F \cup \{f\} \subset \mathcal R_1.$  Pick a finite set of selfadjoint generators $G$ for $\mathcal R_1$ with $f \in G.$  $F^{\prime \prime} \subset G^{\prime \prime}.$  It is clear that $X\cup G$ generates $M.$  By [11] and the hyperfinite inequality for $\delta_0$ 

\begin{eqnarray*} \delta_0(X \cup G ) & \leq & \delta_0(X \cup G \cup F \cup \{f\}) \\ & \leq & \delta_0(X \cup F \cup \{f\}) + \delta_0(G \cup F \cup \{f\}) - \delta_0(F \cup \{f\}) \\ & \leq & \delta_0(X \cup F \cup \{f\}) + 1-1 \\ & \leq & \delta_0(X \cup F \cup \{f\}). \\
\end{eqnarray*}

\noindent Thus, $\delta_0(X \cup G) \leq \delta_0(X \cup F \cup \{f\}) \leq \delta_0(X \cup F).$
\end{proof}

\section{Amplifications and Compressions}

In this section we will compute formulae or bounds relating the free entropy dimension of a set of generators for a $II_1$-factor and corresponding sets of generators for an amplification or cutdown of the factor.  We start with the rational case and then use this to derive weaker results for the general case.  

\subsection{The Rational Case}

Throughout this section denote by $x_1, \ldots, x_n$ a finite set of elements in a tracial von Neumann algebra $A$ and by $\langle e_{ij} \rangle_{1 \leq i,j \leq p}$ the canonical system of matrix units for $M_p(\mathbb C).$  Set $X = \{x_1 \otimes I, \ldots, x_n \otimes I\} \subset A \otimes M_p(\mathbb C)$ and $Y = \langle I \otimes e_{ij} \rangle_{1\leq i,j \leq p} \subset A \otimes M_p(\mathbb C).$  Denote by $\psi$ the canonical tracial state on $A \otimes M_p(\mathbb C).$  Observe that $\delta_0(X) = \delta_0(x_1,\ldots, x_n)$ and $\delta_0(Y) = \delta_0(\langle e_{ij}\rangle_{1 \leq i,j\leq p}).$ Our main goal in this section is to show:

\begin{lemma}$\delta_0( X \cup Y) = 1 - \frac{1}{p^2} + \frac{1}{p^2} \cdot \delta_0(x_1,\ldots, x_n).$ \end{lemma}

\begin{proof} Given $k \in \mathbb N$ denote by $\pi_k: M_p(\mathbb C) \rightarrow M_k(\mathbb C)$ the map which sends an arbitrary $x \in M_p(\mathbb C)$ to the $k \times k$ matrix obtained by taking each entry of $x_{ij}$ and stretching it into the $[ \frac{k}{p}] \times [ \frac{k}{p}]$ matrix with $x_{ij}$ repeated along the diagonals and $0$'s elsewhere.  It is clear that $\pi_k$ is a $*$-homomorphism.  Set $\xi_k = \langle \pi_k(e_{ij})\rangle_{1 \leq i,j \leq p}.$  Given $m \in \mathbb N$ and $\gamma > 0,$ $\xi_k \in \Gamma(Y;m,k,\gamma)$ for $k$ large enough.

Recall the the relative microstate spaces defined in [5]: $\Xi(X;m,k,\gamma),$ $K_{\epsilon}(\Xi(X)),$ and $\mathbb P_{\epsilon}(\Xi(X))$ relative to the sequence $\langle \xi_k \rangle_{k=1}^{\infty}$ given in the preceding paragraph.  By [4] and [5]

\[ \delta_0(X \cup Y) = 1 - \frac{1}{p^2} + \limsup_{\epsilon \rightarrow 0} \frac{\mathbb P_{\epsilon}(\Xi(X))}{|\log \epsilon|} = 1 - \frac{1}{p^2} + \limsup_{\epsilon \rightarrow 0} \frac{\mathbb K_{\epsilon}(\Xi(X))}{|\log \epsilon|}.  \]

\noindent Hence, it suffices to show that $\limsup_{\epsilon \rightarrow 0} \frac{\mathbb K_{\epsilon}(\Xi(X))}{|\log \epsilon|} \leq \frac{1}{p^2} \cdot \delta_0(X)$ and $\limsup_{\epsilon \rightarrow 0} \frac{\mathbb P_{\epsilon}(\Xi(X))}{|\log \epsilon|} \geq \frac{1}{p^2} \cdot \delta_0(X).$

Concerning the first claim given $t>0$ there exists an $\epsilon_0 >0$ such that for $\epsilon_0 > \epsilon >0,$ $\mathbb K_{\epsilon}(X) < (\delta_0(x_1,\ldots, x_n) + t) |\log \epsilon|.$ Find $m > 10$ and $\gamma < (p+1)^{-3}\epsilon^2 $ such that for $k$ sufficiently large there exists an $\epsilon$-net  
$\langle (z^{(1)}_{rk},\ldots, z^{(n)}_{rk}) \rangle_{r \in \Omega(k)}$ for $\Gamma(x_1,\ldots, x_n;m,k,\gamma)$ with $k^{-2} \cdot \log \left ( \# \Omega(k) \right) \leq |\log \epsilon| \cdot (\delta_0(x_1,\ldots, x_n) + 2t).$  Denote by $b^{(l)}_{rk}$ the $k \times k$ matrix which has $z^{(l)}_{r[\frac{k}{p}]}$ repeated along the diagonal $p$ times and $0$ in all other entries.
I claim that $\langle (b^{(1)}_{rk},\ldots, b^{(n)}_{rk}) \rangle_{r  \in \Omega([\frac{k}{n}])}$ is an $4 p^3 \epsilon$-cover for $\Xi(X;m,k,2\gamma).$

Towards this end suppose $a = (a_1,\ldots, a_n) \in \Xi(X;m,k,\gamma).$  For any $*$-monomial f in $n$ noncommuting variables 
\[ \psi \left (f(x_1 \otimes e_{ii}, \ldots, x_n \otimes e_{ii} )\right) = \frac{1}{p} \cdot \psi (\left(f(x_1 \otimes I, \ldots, x_n \otimes I) \right).\]

\noindent Thus,
\begin{eqnarray} p \cdot (\pi_k(e_{11})a_1 \pi_k(e_{11}), \ldots, \pi_k(e_{11})a_n\pi_k(e_{11})) \in \Gamma(x_1,\ldots, x_n;m,k, 2 \gamma). \end{eqnarray}

\noindent For any $1 \leq i \leq p$ and $1 \leq l \leq n,$  $(I \otimes e_{ij})(x_l \otimes I) (I \otimes e_{ji}) = x_1 \otimes e_{ii}$ so that  

\begin{eqnarray} |\pi_k(e_{i1}) a_l \pi_k(e_{1i}) -  a_l \pi_k(e_{ii})|_2 < \epsilon.
\end{eqnarray}

\noindent Finally, for any $1 \leq l \leq n,$ $x_l - \sum_{i=1}^p x_l (I \otimes e_{ii}) = 0$ so that 

\begin{eqnarray} |a_l - \sum_{i=1}^p a_l \pi_k(e_{ii})|_2 < \epsilon. 
\end{eqnarray}

For each $1 \leq l \leq n$ denote by $\eta_l$ the $[\frac{k}{p}] \times [\frac{k}{p}]$ submatrix of $\pi_k(e_{11})a_l\pi_k(e_{11}),$ i.e., the matrix formed from the nonzero entries of $\pi_k(e_{11})a_l\pi_k(e_{11}).$  It is clear from (1) that $\eta = (\eta_1,\ldots, \eta_n) \in \Gamma(x_1,\ldots, x_n;m,[\frac{k}{p}],\gamma).$  There exists some $r \in \Omega([\frac{k}{p}])$ such that $|\zeta_{rk} - \eta|_2 < \epsilon.$  Using (1) and (2) we get

\begin{eqnarray*} |a_l - b^{(l)}_{rk}|_2 \leq |\sum_{i=1}^p a_l \pi_k(e_{ii}) - b^{(l)}_{rk}|_2 + \epsilon & = & |\sum_{i=1}^p  \left(a_l \pi_k(e_{ii}) - b^{(l)}_{rk} \pi_k(e_{ii})\right)|_2 + \epsilon \\ & \leq & \sum_{i=1}^p |\pi_k(e_{i1})a_l\pi_k(e_{1i}) - \pi_k(e_{ii}) b^{(l)}_{rk} \pi_k(e_{ii})|_2 + (p+1) \epsilon \\ & \leq & \sum_{i=1}^p |\pi_k(e_{11}) a_l \pi_k(e_{11}) - \pi_k (e_{11}) b^{(l)}_{rk} \pi_k(e_{11})|_2 + 2(p+1) \epsilon \\ & \leq & \sum_{i=1}^p p \cdot |\eta_l - z_{lk}|_2 + 3(p+1)  \\ & \leq & 4p^2 \epsilon.
\end{eqnarray*}

\noindent It follows that $|a - (b^{(1)}_{rk},\ldots, b^{(n)}_{rk})|_2 < 4p^3 \epsilon$ as desired.  Thus, 

\[ \mathbb K_{4p^3\epsilon}(\Xi(X)) \leq \mathbb K_{4p^3\epsilon}(\Xi(X;m,2\gamma)) \leq \limsup_{k \rightarrow \infty} k^{-2} \cdot \log (\# \Omega_{[\frac{k}{p}]}) \leq \frac{|\log \epsilon|}{p^2} \cdot (\delta_0(x_1,\ldots, x_n) + 2t). \]

\noindent Hence, for $\epsilon$ sufficiently small,

\[ \frac {\mathbb K_{\epsilon}(\Xi(X))}{|\log \epsilon|} \leq \frac{1}{p^2} \cdot \delta_0(x_1,\ldots, x_n) + 2t, \]

\noindent and $t >0$ being arbitrary we have that $\limsup_{\epsilon \rightarrow 0} \frac{\mathbb K_{\epsilon}(\Xi(X))}{|\log \epsilon|} \leq \frac{1}{p^2} \cdot \delta_0(x_1,\ldots, x_n).$  This completes the upper bound claim.

For the lower bound claim suppose $t >0.$  For any $\epsilon_0 >0$ there exists an $\epsilon >0$ such that $\epsilon_0 > \epsilon$ and $\mathbb P_{\epsilon}(X) \geq |\log \epsilon| \cdot (\delta_0(x_1,\ldots,x_n) + t).$  Suppose now that $m \in \mathbb N$ and $\gamma >0$ are given. For $k > \frac{p}{\gamma}$ sufficiently large there exists an $\epsilon$-separated subset $\langle (y^{(1)}_{sk},\ldots, y^{(n)}_{sk})\rangle_{s \in \Pi_k}$ of $\Gamma(X;m,[\frac{k}{p}], \gamma)$ satisfying $k^{-2} \cdot \log (\# \Pi_k) \geq \frac{1}{p^2} \cdot|\log \epsilon| \cdot \delta_0(x_1,\ldots, x_n).$  For each such $k$ and $1 \leq l \leq n$ denote by $b_{sk}^{(l)}$
the $k \times k$ matrix obtained by repeating $y_{sk}^{(l)}$ repeated $p$ times along the diagonal and $0$'s elsewhere.  Now if $1 \leq j \leq m,$ $1 \leq i_1,\ldots, i_j \leq n,$ and $v_1,\ldots, v_{j+1} \in \xi_k \cup \{I_k\}$ where $I_k$ is the $k \times k$ identity matrix, then

\begin{eqnarray*} tr_k(v_1 b_{sk}^{(i_1)} v_2 b_{sk}^{(i_2)} \cdots v_k b^{(i_j)}_{sk} v_{(k+1)}) & = & tr_k(b_{sk}^{(i_1)} \cdots b_{sk}^{(i_j)} v_1 \cdots v_{k+1}) \\ & = & p \cdot tr_{[\frac{k}{p}]}(b_{sk}^{(i_1)} \cdots b_{sk}^{(i_j)}) \cdot tr_k(v_1 \cdots v_{k+1}) 
\end{eqnarray*}

\noindent From this it follows that $\langle (b_{sk}^{(1)}, \ldots, b_{sk}^{(l)}) \rangle_{s \in \Pi_k}$ is an $\epsilon$-separated subset of $\Xi(X;m,k,p\gamma).$  So

\begin{eqnarray*} \mathbb P_{\epsilon}(\Xi(X;m,p\gamma)) & = & \limsup_{k \rightarrow \infty} k^{-2} \cdot \log \left (P_{\epsilon}(\Xi(X;m,k,p\gamma)) \right) \\ & \geq & \limsup_{k \rightarrow \infty} k^{-2} \cdot \log (\# \Pi_k) \\ & \geq & \frac{1}{p^2} \cdot |\log \epsilon| \cdot (\delta_0(x_1,\ldots, x_n) +t)
\end{eqnarray*}

\noindent This being true for any $m \in \mathbb N,$ and $\gamma >0,$ $\frac{\mathbb P_{\epsilon}(\Xi(X))}{|\log \epsilon|} \geq \frac{1}{p^2} \cdot (\delta_0(x_1,\ldots, x_n) +t).$  It follows that $\limsup_{\epsilon \rightarrow 0} \frac{\mathbb P_{\epsilon}(\Xi(X))}{|\log \epsilon|} \geq \frac{1}{p^2} \cdot \delta_0(x_1,\ldots,x_n).$  This completes the lower bound claim.
\end{proof}

\begin{corollary} Suppose $e \in A$ is a projection with trace value $\frac{1}{p}$ for some $p \in \mathbb N,$ $Y = \langle f_{ij} \rangle_{1 \leq i, j \leq n}$ is a system of matrix units in $A$ with $f_{11} = e$, and $X$ is a finite set of generators for $eAe.$  Then $Z = X \cup Y$ is a finite set of generators for $A$ such that 

\[\delta_0(X) = p^2(\delta_0(Z) - 1) +1\]

\noindent where $\delta_0(X)$ is computed in the von Neumann algebra $eAe.$
\end{corollary} 
\begin{proof} Clearly $Z$ generates $A.$  Also, $A \simeq eAe \otimes M_p(\mathbb C)$ via a map which sends $x \in eAe \mapsto x \otimes e_{11}$ and $f_{ij} \mapsto I \otimes e_{ij}$ where $\langle e_{ij} \rangle_{1\leq i,j\leq n}$ is the canonical system of matrix units for $M_p(\mathbb C).$  Call the image of $Z$ in $eAe \otimes M_p(\mathbb C)$ $D.$  Then clearly $D$ generates the same algebra as $\{ x \otimes I \in eAe \otimes M_p(\mathbb C): x \in eAe\} \cup \langle I \otimes e_{ij} \rangle_{1\leq i,j \leq n}$ so that by the previous lemma, $\delta_0(Z) = \delta_0(D) = 1 - \frac{1}{p^2} + \frac{1}{p^2} \cdot \delta_0(X).$  
\end{proof}

Recall that for each $\alpha \in (0, \infty)$ there exists a type $II_1$-factor $M_{\alpha}$ ($M$ as always, a $II_1$-factor) obtained by taking a projection $e \in M \otimes B(\ell^2(\mathbb Z))$ satisfying $Tr(e) = \alpha$ ($Tr$ the canonical trace) and setting $M_{\alpha} = eMe.$ 

\begin{corollary} Suppose $\alpha \in \mathbb Q_{>0}.$  If $M$ can be finitely generated then there exist two finite sets of generators $F$ for $M$ and $G$ for $M_{\alpha},$ such that

\[ \delta_0(G) = 1 - \frac{1}{\alpha^2} + \frac{1}{\alpha^2} \cdot \delta_0(F).\] 
\end{corollary}

\begin{proof} There exist natural numbers $m$ and $n$ such that $\alpha = \frac{m}{n}.$  Because $M$ can be generated by finitely many elements, so 
can $M_{\frac{1}{n}}.$  Fix such a set of generators $X$ for $M_{\frac{1}{n}}$ canonically identified in $M_{\frac{1}{n}} \otimes M_m(\mathbb C)$ and denote by $E$ the obvious system of matrix units for $M_m(\mathbb C)$ canonically identified in $M_{\frac{1}{n}} \otimes M_m(\mathbb C).$  Set $G = X \cup E$ and observe that the von Neumann algebra generated by $G$ is $M_{\frac{1}{n}} \otimes M_m(\mathbb C) = M_{\alpha}.$  From Lemma 3.1 we have that

\[ \delta_0(G) = 1 - \frac{1}{m^2} + \frac{1}{m^2} \cdot \delta_0(X).
\]

\noindent By Corollary 3.2 there exists a finite set of generators $F$ for $M$ such that $\delta_0(X) = n^2(\delta_0(F) -1)+1.$  Substituting this into the above equation gives the result.
\end{proof}

\subsection{The General Case} Suppose again that $\infty > \alpha >0.$  As always, $M$ is a $II_1$-factor.

\begin{lemma} Suppose $M$ can be generated by a finite number of elements and $0 < \alpha <\infty.$ For $\epsilon >0$ there exist finite sets of generators $F$ and $G$ for $M$ and $M_{\alpha},$ respectively, such that

\[ \delta_0(G) < 1 - \frac{1}{\alpha^2} + \left (\frac{1}{\alpha^2} + \epsilon \right )\cdot \delta_0(F) + \epsilon.
\] 
\end{lemma}

\begin{proof} $M_{\alpha} = e(M \otimes B(H))e$ for some finite dimensional Hilbert space and $e$ is a projection in $M \otimes B(H)$ with trace $\alpha.$  Suppose $n \in \mathbb N$ is a multiple of $\dim H.$  For some $j \in \mathbb N,$ $\dim H \cdot \frac{j}{n} \leq \alpha \leq \dim H \cdot \frac{j+1}{n}.$  Choose a system of matrix units in $ \langle e_{pq} \rangle_{1 \leq p,q\leq n}$ in $M \otimes B(H).$ It can be arranged so that $e_{11} + \cdots + e_{jj} \leq e$ and that the algebra the matrix units form is $A \otimes B(H)$ where $A$ is a type $I$ finite factor in $M.$ Set $f_n = e - (e_{11} + \cdots + e_{jj}).$  Because $M$ can be generated by finitely many elements, $e_{11}Me_{11}$ can also.  Fix a finite set of generators $X_n$ for $e_{11}Me_{11}.$ In this same way we fix such a finite set of generators $Y$ for $M_{\alpha}.$  Set $E_n = \langle e_{pq} \rangle_{1 \leq p,q \leq j}.$Now consider the finite set

\[ G_n= X_n \cup E \cup Yf_n \cup f_nY \cup \{f_n\}
\]  

\noindent where $Yf_n = \{yf: y \in Y\}$ and $f_nY= \{f_ny: y \in Y\}.$  $G_n$ is a set of generators for $M_{\alpha}.$ 

I will bound $\delta_0(G_n)$ by appealing again to the relative (or conditional) microstate space result in [5].  Although it was phrased in terms of selfadjoint microstate spaces, it holds true in the non-selfadjoint case as well.  We condition with respect to microstates for $f_n.$  For each $k$ denote by $p_k$ the diagonal $k\times k$ matrix with $1$'s on the last $[\varphi(f_n)]k$ entries and $0$'s everywhere else.  Observe that for any $m \in \mathbb N$ and $\gamma >0,$ for sufficiently large $k,$ $p_k \in \Gamma(f_n;m,k,\gamma).$  By [5]

\begin{eqnarray} \delta_0(G_n) = \delta_0(f_n) + \limsup_{\epsilon \rightarrow 0} \frac{\mathbb K_{\epsilon}(\Xi(X_n \cup \langle e_{pq}\rangle_{1 \leq p,q \leq j} \cup Yf_n \cup f_nY))}{|\log \epsilon|}
\end{eqnarray}

\noindent where the covering quantity is taken relative to $\langle p_k \rangle_{k=1}^{\infty}.$

Suppose $t >0.$  By Corollary 3.2 $\delta_0(X_n \cup E_n) = 1 - \frac{1}{j^2} + \frac{1}{j^2} \cdot \delta_0(X_n).$  Here  $\delta_0(X)$ is computed in the von Neumann algebra $e_{11}Me_{11}$ and $\delta_0(X_n \cup E_n)$ is computed in the von Neumann algebra $e_{11}Me_{11} \otimes M_n(\mathbb C).$   Thus, for $\epsilon >0$ sufficiently small there exist $m \in \mathbb N$ and $\gamma >0$ (dependent upon $\epsilon$) and an $\epsilon$-cover of $\Gamma(X_n \cup E_n;m,k,2\gamma)$ such that for $k$ sufficiently large the logarithm of the  cardinality of the cover is greater than

\[  |\log \epsilon| \cdot k^2 \cdot \left(1 - \frac{1}{j^2} + \frac{1}{j^2} \cdot \delta_0(X_n) + t \right). \]

\noindent We can also choose $m$ and $\gamma$ so that if $(a,b)\in \Xi(Yf_n \cup f_nY  ;m,k,\gamma)$ (relative to the sequence $\langle p_k \rangle_{k=1}^{\infty}$), then $|a p_k - a|_2, |p_k b-b|_2 < \epsilon,$ where $a p$ and $p b$ are the tuples obtained by multiplying each entry from the left or right, respectively, by $p.$  Finally, there exists $m_1 > m$ and $\gamma_1 < \gamma$ such that if $a \in \Xi(X_n \cup E_n;m_1,k,\gamma_1)$ (microstate space relative to $\langle p_k \rangle_{k=1}^{\infty}),$ then $p_k^{\bot} a p_k^{\bot}$ (the $(\#X_n + j^2)$-tuple obtained by compressing each entry by $p_k$) is in $\Xi(X_n \cup E;m,k,2\gamma)$ and $|p_k^{\bot} a p_k^{\bot} - a|_2 < \epsilon.$ 

Now for each $k$ find an $\epsilon$-cover $\langle \xi_{rk} \rangle_{r \in \Omega_k}$ of $\Gamma(X_n \cup E_n;m, k-[\varphi(f)k], \gamma).$  For each such $r$ and $k$ denote by $a_{rk}$ the $(\#X_n +j^2)$-tuple whose $i^{th}$ entry is the $k\times k$ matrix which the $i^{th}$ component of $\xi_{rk}$ as a submatrix in the left corner and $0$'s elsewhere.   For each $k$ find an $\epsilon$-net $\langle v_{sk} \rangle_{s \in \Pi_{k}}$ in the $|\cdot|_2$-ball of $(M_k(\mathbb C))^{\#Y}p_{k} \bigoplus p_k (M_k(\mathbb C))^{\#Y},$ of radius $K= |(Y_n,Y_n)|_2 +1.$  $(M_k(\mathbb C))^{\#Y}p_k \bigoplus p_k(M_k(\mathbb C))^{\#Y} \subset (M_k(\mathbb C))^{2\#Y}p_k$ being Euclidean space of dimension no greater than $4 k^2 \varphi(f) ,$ the standard volume comparison test allows us to assume that $\# \Pi_k < \left (\frac{K+1}{\epsilon} \right)^{4k^2\varphi(f) }.$  

I claim that $\langle (a_{rk}, v_{sk}) \rangle_{(r,s) \in \Omega_k \times \Pi_k}$ is a $4\epsilon$-cover for $\Xi(X_n \cup E \cup fY_n \cup Yf_n;m_1,k,\gamma_1).$  To see this suppose $(x,y,z) \in \Xi(X_n \cup E_n \cup Yf_n \cup f_nY;m,k,\gamma)$ where $x$ is a $(\#X_n + j^2)$-tuple (corresponding to $X_n \cup E$), $y$ is a $\#Y_n$-tuple corresponding to $Yf,$ and $z$ is  $\#Y_n$-tuple corresponding to $fY_n$.  By the choice of $m_1$ and $\gamma_1,$ 
$|x - p_k^{\bot}x p_k^{\bot}|_2 < \epsilon$ and $p_k^{\bot}xp_k^{\bot} \in \Xi(X_n \cup E_n;m,k,2\gamma).$  It follows from the previous paragraph that there exists an $r_0 \in \Omega_k$ such that $|a_{r_0k} - p_k^{\bot}x p_k^{\bot}|_2 < \epsilon,$ whence $|a_{r_0k} - x|_2 < 2 \epsilon.$  Again by the choice of $m_1$ and $\gamma_1,$ $|yp_k -y|_2, |p_kz -z|_2 <\epsilon.$  There exists an $s_0 \in \Pi_k$ satisfying $|(yp_k, p_kz) - v_{s_0k}|_2 < \epsilon$ whence $|(y,z) - v_{s_0k}|_2 < 3 \epsilon.$  Thus, $|(x,y,z) - (a_{r_0k}, v_{s_0k})|_2 < 4\epsilon.$
From the preceding paragraph we now have:

\begin{eqnarray*} \mathbb K_{5 \epsilon}(\Xi(X_n \cup E_n \cup Yf_n \cup f_nY)) & \leq &  \mathbb K_{5 \epsilon}(\Xi(X_n \cup E_n \cup Yf_n \cup f_nY;m_1,\gamma_1)) \\ & \leq & \limsup_{k \rightarrow \infty} k^{-2} \cdot \log(\#\Omega_k \cdot \#\Pi_k) \\ & \leq & \limsup_{k \rightarrow \infty} k^{-2} \cdot \log(\#\Omega_k) + \limsup_{k \rightarrow \infty} k^{-2} \cdot \log (\#\Pi_k) \\ & \leq & (1 - \varphi(f))^2 \cdot |\log \epsilon| \cdot \left(1 - \frac{1}{j^2} + \frac{1}{j^2} \cdot \delta_0(X) + t \right) + \\ & & |\log \epsilon| \cdot 4 \varphi(f_n)
\end{eqnarray*}

\noindent This being true for all $\epsilon$ sufficiently small and $t >0$ being arbitrary (4) yields 

\[ \delta_0(G_n) \leq \delta_0(f_n) + (1 -\varphi(f_n))^2 \cdot \left (1 -\frac{1}{j^2} + \frac{1}{j^2}\cdot \delta_0(X_n) \right) + \varphi(f_n).\]

By Corollary 3.2 we can find a finite set of generators $F_n$ for $M$ such that 

\[\delta_0(X_n) = \left (\frac{n}{\dim H} \right )^2 (\delta_0(F_n) -1) +1.\]

\noindent Substituting this above gives

\[ \delta_0(G_n) \leq \delta_0(f_n) + 1 + \frac{n^2}{(j \cdot \dim H)^2} \cdot (\delta_0(F_n) - 1) + \varphi(f_n)
\]

\noindent As $n \rightarrow \infty,$ $\varphi(f_n) \rightarrow 0,$ $\delta_0(f_n) \rightarrow 0,$ and $\frac{n^2}{(j \cdot \dim H)^2} \rightarrow \frac{1}{\alpha^2}.$  This yields the desired result.  
\end{proof}

\begin{corollary} If $M$ can be generated by finitely many elements and $0 < \alpha < \infty,$ then for $\epsilon >0$ there exist finite sets of generators $F$ and $G$ for $M$ and $M_{\alpha},$ respectively such that

\[ \delta_0(F) < 1 - \alpha^2 + (\alpha + \epsilon)^2 \cdot \delta_0(G) + \epsilon.
\]
\end{corollary}

\section{Weak Bounds for an inclusion of $II_1$-factors of finite index}

In this section I will use the results of the previous two sections to derive the main inequality concerning $\delta_0$ and finite index subfactors.  $N \subset M$ is an inclusion of $II_1$-subfactors with $[M:N] < \infty.$  Moreover, assume $N$ can be generated by finitely many elements.  First for the rational case:

\begin{lemma} Suppose $[M:N] \in \mathbb Q.$  Then there exist finite generating sets $X$ for $N$ and $Y$ for $M$ such that

\[ \delta_0(X) \geq \delta_0(Y) \geq [M:N]^{-2} \cdot (\delta_0(X)-1) +1.
\]
\end{lemma}

\begin{proof} Write $[M:N] = \frac{m}{n}$ where $m, n \in \mathbb N.$  Find a copy $\mathcal R$ of the hyperfinite $II_1$-factor in $N_1$ and a system of matrix units $E = \langle e_{ij} \rangle_{1 \leq i,j \leq n}$ in $\mathcal R.$  Choose a finite set of selfadjoint generators $D$ for $\mathcal R.$ Because $N$ can be finitely generated so can $e_{11}Ne_{11}.$  Choose a finite set of generators $Z$ for $e_{11}Ne_{11}.$  Set $X = Z \cup E \cup D.$  $X$ clearly generates $N.$  Lemma 2.3 implies

\[ \delta_0(X) \geq \delta_0(X \cup G)\]

\noindent for some finite set of selfadjoint elements $G$ satisfying $f_0 \in G,$ $ \mathcal R \subset G^{\prime \prime} \simeq \mathcal R$ , and $(X \cup G)^{\prime \prime}=M.$ 

Denote by $e_N \in B(L^2(M))$ the projection onto $L^2(N) \subset L^2(M).$ Lemma 2.3 applied to the inclusion $M \subset M_1$ and the generating set $X \cup G$ for $M$ yields a finite set of selfadjoint elements $W$ such that $\mathcal R \subset W^{\prime \prime} \simeq R,$ $e_N \in W,$ and $\delta_0(X\cup G) \geq \delta_0(X \cup G \cup W ).$  

Because $e_N \in N^{\prime},$ $e_{11}e_N$ is a projection of trace $\frac{1}{m}$ in the hyperfinite $II_1$-factor $W^{\prime \prime}.$  Find a system of matrix units $\langle f_{ij} \rangle_{1 \leq i,j \leq m}$ in $W^{\prime \prime}$ such that $f_{11} =e_{11}e_N.$ Clearly

\[ f_{11}M_1f_{11} =e_{11}e_N M_1 e_Ne_{11} = e_{11}Ne_{11}e_N \simeq e_{11}Ne_{11}.\]
 
\noindent It follows from Corollary 3.2 that $H = e_NZe_N \cup \langle f_{ij} \rangle_{1 \leq i,j \leq m}$ is a generating set for $M_1$ such that $\delta_0(Z) = m^2 (\delta_0(H) -1) +1$ where $\delta_0(Z)$ is computed in the von Neumann algebra $e_{11}Ne_{11}.$  By Corollary 4.2 of [5], algebraic invariance of $\delta_0$ ([11]), and Corollary 3.2

\begin{eqnarray*} \delta_0(X \cup G \cup W) & = & \delta_0(Z \cup E \cup G \cup W) \\ & = & \delta_0(Z \cup E \cup G \cup W \cup \langle f_{ij} \rangle_{1 \leq i,j \leq m}) \\ & = & \delta_0(Z \cup Ze_N \cup E \cup G \cup W \cup \langle f_{ij} \rangle_{1 \leq i,j \leq m}) \\ & \geq & \delta_0(e_NZe_N \cup \langle f_{ij} \rangle_{1 \leq i,j \leq m}) \\ & = & \frac{1}{m^2} \cdot \delta_0(Z)- \frac{1}{m^2} +1\\ & = & [M:N]^{-2} \cdot (\delta_0(X) - 1) + 1.\\
\end{eqnarray*}

   Set $Y = X \cup G.$  $Y$ generates $M$ and by the first and second paragraphs and the computations above

\[\delta_0(X) \geq \delta_0(Y) = \delta_0(X \cup G) \geq \delta_0(X \cup G \cup W) \geq [M:N]^{-2} \cdot (\delta_0(X) -1) +1.\]

\end{proof}  

We can drop the condition that $[M:N] \in \mathbb Q$ in the preceding lemma at the cost of obtaining an approximate version of the resultant inequality: 

\begin{lemma} Suppose $[M:N] < \infty.$  For each $n \in \mathbb N$ then there exist finite generating sets $X_n$ for $N$ and $Y_n$ for $M$ and  sequences of real numbers $\langle \alpha_n \rangle_{n=1}^{\infty}$  and $\langle \beta_n \rangle_{n=1}^{\infty}$ with $\alpha_n \rightarrow [M:N]^{-1}$ and $\beta_n \rightarrow 0$ such that

\[ \delta_0(X_n) \geq \delta_0(Y_n) \geq \alpha_n^2 \cdot (\delta_0(X_n)-1) +1 - \beta_n.
\]
\end{lemma}

\begin{proof} Write $[M:N] = \alpha$ and assume $\epsilon >0$ is given.  Suppose $\alpha < n \in \mathbb N.$  Find the largest $m \in \mathbb N$ satisfying $\frac{m}{n} \leq \alpha^{-1} < \frac{m+1}{n}$ and set $\alpha_n = \frac{m}{n}.$ Consider $\mathcal R \subset N_1$ in Section 2 and fix a finite set of generators $F$ for $\mathcal R$ ($F$ is independent of $n$).  Find a system of matrix units $E_n = \langle e_{ij} \rangle_{1 \leq i,j \leq m}$ in $\mathcal R$ with $\varphi(e_{11}) = \frac{\alpha}{n}$ (the matrix algebra generated by $E_n$ will not contain the identity of $N$).  Because $N$ can be finitely generated so can $e_{11}Ne_{11}.$  Set $f_n = I - e_{11} + \cdots + e_{mm}.$  Choose a finite set of generators $Y$ for $N$ and $Z_n$ for $e_{11}Ne_{11}\subset N$ and arrange it so that $F \subset Y$ and $e_{11}Fe_{11} \subset Z_n.$  Define $X_n = Z_n \cup E_n \cup f_nY \cup Yf_n \cup \{f_n\}.$  If $n$ is large enough, then the proof of Lemma 3.4 implies that 

\[\delta_0(X_n) \leq 1 - \frac{1}{m^2} + \frac{1}{m^2} \cdot \delta_0(Z_n) + \epsilon.\]  

\noindent Whence, $\delta_0(Z_n) \geq m^2 \cdot \delta_0(X_n) - m^2 + 1 - m^2 \epsilon.$

Lemma 2.3 implies

\[ \delta_0(X_n) \geq \delta_0(X_n \cup G_n)\]

\noindent for some finite set of selfadjoint elements $G_n$ satisfying $f \in G_n,$ $ \mathcal R \subset G_n^{\prime \prime} \simeq \mathcal R$ , and $(X \cup G_n)^{\prime \prime}=M.$  Here $f$ is the projection onto $L^2(N_1) \subset L^2(N)$ where $N_1 \subset N$ is as in the beginning of section 2. 

Denote by $e_N \in B(L^2(M))$ the projection onto $L^2(N) \subset L^2(M).$ Lemma 2.3 applied to the inclusion $M \subset M_1$ and the generating set $X \cup G_n$ for $M$ yields a finite set of selfadjoint elements $W_n$ such that $\mathcal R \subset W_n^{\prime \prime} \simeq R,$ $e_N \in W,$ and $\delta_0(X_n \cup G_n ) \geq \delta_0(X_n \cup G_n \cup W_n ).$  

Because $e_N \in N^{\prime},$ $e_{11}e_N$ is a projection of trace $\frac{1}{n}$ in the hyperfinite $II_1$-factor $W_n^{\prime \prime}.$  Find a system of matrix units $\langle f_{ij} \rangle_{1 \leq i,j \leq n}$ in $W_n^{\prime \prime}$ such that $f_{11} = e_{11}e_N.$  Clearly

\[ f_{11}M_1f_{11} =e_{11}e_N M_1 e_Ne_{11} = e_{11}Ne_{11}e_N \simeq e_{11}Ne_{11}.\]
 
\noindent It follows from Corollary 3.2 that $H_n = e_N Z_n e_N \cup \langle f_{ij} \rangle_{1\leq i,j \leq n}$ is a finite set of generators for $M_1$ such that $\delta_0(Z_n) = n^2 (\delta_0(H_n) -1) +1$
where $\delta_0(Z_n)$ is computed in the von Neumann algebra $e_{11}Ne_{11}.$  By Corollary 4.2 of [5], algebraic invariance of $\delta_0$ ([11]), and Corollary 3.2

\begin{eqnarray*} \delta_0(X_n \cup G_n \cup W_n) & = & \delta_0(X_n \cup G_n \cup W_n \cup \langle f_{ij} \rangle_{1 \leq i,j \leq n}) \\ & = &  \delta_0(X_n \cup Z_n e_N \cup G_n \cup W_n \cup \langle f_{ij} \rangle_{1 \leq i,j \leq n}) \\ & \geq & \delta_0(e_NZ_ne_N \cup \langle f_{ij} \rangle_{1 \leq i,j \leq n}) \\ & = & \frac{1}{n^2} \cdot \delta_0(Z_n)- \frac{1}{n^2} +1\\ & \geq & \alpha_n^2 \cdot (\delta_0(X_n) -1) + 1 - \alpha_n^2  \cdot \epsilon \\ 
\end{eqnarray*}

   Set $Y_n = X_n \cup G_n.$  $Y_n$ generates $M$ and by what preceded

\[\delta_0(X_n) \geq \delta_0(Y_n)\geq \delta_0(X_n \cup G_n \cup W_n) \geq \alpha_n^2 \cdot (\delta_0(X_n) -1) +1 - \alpha_n^2 \cdot \epsilon.\]

\noindent Since $\epsilon >0$ was arbitrary we are done. 

\end{proof}  

\section{Applications}

We conclude with a few simple applications of the preceding section.
Suppose $L(F_r)$ is the interpolated free group factor with index $r>1$ ([1], [7]).

\begin{corollary} If $L(F_r)$ is a finite index subfactor of $M,$ then there exists a finite generating set $Y$ for $M$ such that $\delta_0(Y) \geq [M:L(F_r)]^{-2} \cdot (r-1) + 1.$ 
\end{corollary}

\begin{proof}  Set $N = L(F_r)$ and recall from the proof of the Lemma 4.2 where $Z_n$ is chosen to be a set of generators for $eNe$ where $e \in N$ is a projection of trace $[M:N]/n.$  By results of Radulescu and Dykema ([1], [7]), $e  L(F_r) e \simeq L(F_s)$ where $s =\frac{n^2}{[M:N]^2} \cdot (r-1) +1$ and we may choose $Z_n$ so that $\delta_0(Z_n) \geq \frac{n^2}{[M:N]^2} \cdot (r-1) +1.$  Constructing the generating set $Y=Y_n$ for $M$ with respect to this $Z_n$ we have again by the proof of Lemma 4.2 that

\[ \delta_0(Y) \geq \frac{1}{n^2} \cdot \delta_0(Z_n) - \frac{1}{n^2} +1 
 \geq [M:L(F_r)]^{-2} \cdot (r-1) + 1.
\]

\end{proof}

Using the basic construction the same result can be obtained concerning finite index subfactors of $L(F_r).$  First for a preliminary lemma:

\begin{lemma} Suppose there exists a $t >0$ such that for any finite set of generators $F$ for $M,$ $\delta_0(F) \leq t.$  If $\alpha >0,$ then for any finite set of generators $G$ for $M_{\alpha},$ $\delta_0(G) \leq 1 - \alpha^{-2} + t \alpha^{-2}.$
\end{lemma}

\begin{proof}  First consider the case where $\alpha = p^{-1},$ $p \in \mathbb N.$  If $X$ is finite set of generators for $M_{p^{-1}}$ then consider the set of generators $Y$ for $M_{p^{-1}} \otimes M_p(\mathbb C) =M$ consisting of elements of the form $I \otimes e_{ij},$ $1 \leq i,j \leq p$ where the $e_{ij}$ are the canonical system of matrix units for $M_p(\mathbb C)$ and elements of the form $x \otimes e_{11}$ where $x \in X.$  It follows from Corollary 3.2 and the hypothesis that $1 - p^{-2} + p^{-2} \cdot \delta_0(X) = \delta_0(Y) < t.$  

Now suppose $\alpha$ is arbitrary, $X$ is a set of generators for $M_{\alpha},$ and let $\epsilon >0.$  Suppose $n \in \mathbb N$ and find a system of matrix units $E = \langle e_{ij} \rangle_{1 \leq i,j \leq m}$ such that $\varphi(e_{11}) = (n\alpha)^{-1}$ and $m$ is the largest integer such that $m \cdot (n \alpha)^{-1} \leq 1.$  Notice that this system of matrix units will not contain the identity.  By [11], $\delta_0(X) \leq \delta_0(X \cup E).$  Set $f = I - (e_{11} + \cdots + e_{nn}).$  Denote by $Y$ the set consisting of all elements of the form $e_{1i} x e_{j1}$ where $x \in X.$  $Y$ is a generating set for $e_{11} M e_{11} \simeq M_{n^{-1}}.$  $X \cup E$ and $Y \cup E \cup fX \cup Xf$ generate the same unital $*$-algebra whence by [11] their free entropy dimensions are equal. Set $Z = Y \cup E \cup fX \cup Xf \cup \{f\}.$  Now by the proof of Lemma 3.4 if $n$ is chosen large enough, then

\begin{eqnarray*} \delta_0(X) \leq \delta_0(X \cup E) & \leq & \delta_0(Y \cup E \cup fX \cup Xf) \\ & \leq & \delta_0(Z) \\ & \leq &  1 - \frac{1}{m^2} + \frac{1}{m^2} \cdot \delta_0(Y) + \epsilon \\ & \leq & 1 - \frac{1}{m^2} + \frac{n^2}{m^2} \cdot (t-1) + \frac{1}{m^2} + \epsilon.
\end{eqnarray*}

\noindent As $n \rightarrow \infty$, $m \rightarrow \infty$ and $\frac{n}{m} \rightarrow \alpha^{-1}.$  Thus, $\delta_0(X) \leq 1 + \alpha^{-2}(t-1)+ \epsilon.$  $\epsilon >0$ being arbitrary it follows that $\delta_0(X) \leq 1 - \alpha^{-2} + t\alpha^{-2}.$ 
\end{proof}

\begin{corollary} If $N$ satisfies the property that for some $\epsilon >0$ and any finite set of generators $X$ for $N,$ $\delta_0(X) < r - \epsilon,$ then $L(F_r)$ cannot contain $N$ as a subfactor of finite index.
\end{corollary}

\begin{proof} This follows from duality.  More precisely suppose $L(F_r)$ contains $N$ as a subfactor of finite index.  Perform the basic construction to arrive at the inclusion $N \subset L(F_r) \subset M.$  By [3] $[M:L(F_r)] = [L(F_r):N] = \alpha < \infty.$  On the other hand $M \simeq N_{\alpha}.$  For any finite set of generators $X$ for $N_{\alpha}$ we have by hypothesis and the previous lemma, $\delta_0(X) \leq 1 - \alpha^{-2} + (r-\epsilon) \cdot \alpha^{-2} < \alpha^{-2} (r-1)+1 .$  By Corollary 5.1 $L(F_r)$ cannot be a finite index subfactor of $M$ of index $\alpha.$  This is absurd. 
\end{proof}

Using [2] we get:

\begin{corollary} If $M$ can be generated by a sequence of Haar unitaries $\langle u_j \rangle_{j=1}^{\infty}$ with the property that $u_{j+1}u_j u_{j+1}^* \in \{u_1,\ldots,u_j\}^{\prime \prime}$ for each $j$, then $M$ cannot contain $L(F_r)$ as a finite index subfactor and $L(F_r)$ cannot contain $M$ as a finite index subfactor. 
\end{corollary}

\noindent{\it Acknowledgements.} I thank Dimitri Shlyakhtenko for raising this question about subfactors and free entropy dimension during the 2002 entropy conference at UCLA and subsequent conversations about this problem.  I also would like to thank Dan Voiculescu for some encouraging remarks during the beginning of this work.


\begin{thebibliography}{[ASMR]} 

\bibitem{1} Dykema, Ken {\it Interpolated free group factors.} Pacific Journal of Mathematics, 163 (1994), 123-135.

\bibitem{2} Ge, Liming and Shen, Junhao {\it On free entropy dimension of von Neumann algebras.}, GAFA, (2002)

\bibitem{3} Jones, V.F.R. {\it Index for Subfactors.} Invent. Math. 73 (1983), no.1, 1-25.

\bibitem{4} Jung, Kenley {\it The free entropy dimension of hyperfinite
von Neumann algebras.} Transactions of the AMS, 355 (2003), 5053-5089 .

\bibitem{5} Jung, Kenley {\it A Hyperfinite inequality for $\delta_0$}, preprint.

\bibitem{6} Ozawa, Narutaka {\it Solid von Neumann algebras,} preprint.

\bibitem{7} Radulescu, Florin {\it Random matrices, amalgamated free 
products and subfactors in free group factors of noninteger index.} 
Inventiones mathematicae 115 (1994), 347-389.

\bibitem{8} Shlyakhtenko, Dimitri; Ueda, Yoshimichi {\it Irreducible subfactors of $L(\mathbb F\sb \infty)$ of index $\lambda >4$.} J. Reine Anwgew. Math. 548 (2002), 149 -166 

\bibitem{9} Stefan, Marius {\it The primality of subfactors of finite index in the interpolated free group factors.} Proc. Amer. Math. Soc. 126 (1998), no.8, 2299-2307.

\bibitem{10} Voiculescu, D. {\it The analogues of entropy and of Fisher's
information measure in free probability theory, II}, Inventiones
mathematicae 118, (1994), 411-440.

\bibitem{11} Voiculescu, D. {\it A Strengthened Asymptotic Freeness Result 
for Random Matrices with Applications to Free Entropy}, IMRN, 1 (1998), 
41-64.      

\end{thebibliography}
\end{document}